\documentclass[12pt]{amsart}
\usepackage[latin1]{inputenc}
\usepackage{amsxtra,amssymb,amsthm,amsmath,amscd,mathrsfs, epsfig,eufrak,url}
\usepackage[all]{xy}
\def \Q {{\mathbb Q}}
\def \R {{\mathbb R}}

\def \Z {{\mathbb Z}}

\def \d {\,{\rm d}}

\def\dm{\frac{1}{2}}

\def\eps{\varepsilon}

\def\sumb{\mathop{\sum \Bigl.^{\flat}}\limits}

\newcommand{\stacksum}[2]{\substack{{{#1}}\\{{#2}}}}
\def\le{\leqslant}
\def\leq{\leqslant}
\def\ge{\geqslant}
\def\geq{\geqslant}

\newcommand{\Hf}{\mathrm{H}}

\newcommand{\ra}{\rightarrow}

\newcommand{\mods}[1]{\,(\mathrm{mod}\,{#1})}

\newcommand{\expll}[2]{\exp\Bigl({{#1}}\frac{\log {{#2}}}{\log\log {{#2}}}\Bigr)}

\theoremstyle{plain}
\newtheorem{theorem}{Theorem}

\newtheorem{lemma}{Lemma}[section]
\newtheorem{corollary}{Corollary}

\newtheorem{proposition}{Proposition}

\theoremstyle{remark}
\newtheorem{remark}{Remark}

\theoremstyle{definition}

\numberwithin{equation}{section}

\begin{document}

\title[On modular signs]{On modular signs} 
\author{E. Kowalski, Y.-K. Lau, K. Soundararajan \& J. Wu}
\address{ETH Z\"urich -- D-MATH\\
  R\"amistrasse 101\\
  8092 Z\"urich, Switzerland}
\thanks{The work of E. Kowalski was supported
  in part by the National Science Foundation under agreement No.
  DMS-0635607 during a sabbatical stay at the Institute for Advanced
  Study.}
\email{kowalski@math.ethz.ch}

\address{Department of
  Mathematics, The University of Hong Kong, Pokfulam Road, Hong Kong}
\email{yklau@maths.hku.hk}

\address{Department of Mathematics, Stanford University, Stanford, CA 94305, USA}
\thanks{The work of K. Soundararajan was partially supported by the National 
Science Foundation (DMS 0500711), and through the Veblen Fund 
of the Institute for Advanced Study, Princeton.}
\email{ksound@stanford.edu}

\address{Institut Elie Cartan Nancy (IECN)
  \\
  CNRS, Nancy-Universit\'e, INRIA
  \\
  Boulevard des Aiguillettes, B.P. 239
  \\
  54506 Van\-d\oe uvre-l\`es-Nancy
  \\
  France} \email{wujie@iecn.u-nancy.fr}

\date{\today}

\begin{abstract}
  We consider some questions related to the signs of Hecke eigenvalues
  or Fourier coefficients of classical modular forms.  One problem is
  to determine to what extent those signs, for suitable sets of
  primes, determine uniquely the modular form, and we give both
  individual and statistical results. The second problem, which has
  been considered by a number of authors, is to determine the size, in
  terms of the conductor and weight, of the first sign-change of Hecke
  eigenvalues.  Here we improve the recent estimate of Iwaniec, Kohnen
  and Sengupta.
\end{abstract}



\maketitle


\section{Introduction}

There are many results in the arithmetic of modular forms which are
concerned with various ways of characterizing a given primitive cusp
form $f$ from its siblings, starting from the fact that Fourier
coefficients, hence the $L$-function, determine uniquely a cusp form
$f$ relative to a congruence subgroup $\Gamma$ of $SL(2, \Z)$. Among
such results are stronger forms of the multiplicity one theorem for
automorphic forms or representations, various explicit forms of these
statements, where only finitely many coefficients are required (say at
primes $p\le X$, for some explicit $X$ depending on the parameters
defining $f$), and a number of interesting ``statistic'' versions of
the last problem, where $X$ can be reduced drastically, provided one
accepts some possible exceptions. Among other papers, we can
cite~\cite{GH93}, \cite{M97}, \cite{DK00} or \cite{K05}.
\par
Some of these statements were strongly suggested by the analogy with
the problem of the least quadratic non-residue, which is a problem of
great historic importance in analytic number theory, and there are
many parallels between the results which have been obtained. However,
this parallel breaks down sometimes. For instance, in~\cite{LW08b},
Lau and Wu note that one result of~\cite{LW08a} for the least
quadratic non-residue is highly unlikely to have a good analogue for
modular forms. This result (see~\cite[Th. 3]{LW08a}) is a precise
estimate for the number of primitive real Dirichlet characters of
modulus $q\leq D$ for which the least $n$ with $\chi(n)=-1$ is $\gg
\log D$, and the difficulty is that this estimate can be understood by
assuming that the values $\chi(p)$, for $p$ of moderate size compared
with $D$, behave like independent random variables taking values $\pm
1$ equally often. However, Hecke eigenvalues may take many more than
two values, and thus assuming that they coincide should definitely be
a much more stringent condition.
\par
In this paper, we consider a way to potentially recover a closer
analogy: namely (narrowing our attention to forms with real
eigenvalues) instead of looking at the values of the Hecke
eigenvalues, we consider only their \emph{signs} (where we view $0$ as
being of both signs simultaneously, to increase the possibility of
having same sign). Then classical questions for Dirichlet characters
and modular forms have the following analogues for signs of Hecke
eigenvalues $\lambda_f(p)$ of a classical modular form $f$:
\begin{itemize}
\item What is the first sign-change, i.e., the smallest $n\geq 1$ (or
  prime $p$) for which $\lambda_f(n)<0$ (or $\lambda_f(p)<0$)?
  (Analogues of the least quadratic non-residue). Note a small
  difference with quadratic characters: it is not true here that the
  smallest integer with negative Hecke eigenvalue is necessarily
  prime; finding one or the other are two different
  questions.\footnote{\ E.g., as a random example, for the cusp form
    of weight $2$ associated to the elliptic curve $y^2=x^3+x$, the
    first negative coefficient is $\lambda(9)=-3$, and the first
    negative prime coefficient is $\lambda(13)=-6$.}
\par
\item Given arbitrary signs $\eps_p\in \{\pm 1\}$ for all primes, what
  is the number of $f$ (in a suitable family) for which $\lambda_f(p)$
  has sign $\eps_p$ for all $p\leq X$, for various values of $X$?
  (Analogue of the question in~\cite{LW08a}). 
\par
\item In particular, is there a finite limit $X$ such that
  coincidences of signs of $\lambda_f(p)$ and $\eps_p$ for all $p\leq
  X$ implies that $f$ is uniquely determined? (Analogue of the
  multiplicity one theorem).
\end{itemize}

Of these three problems, only the first one seems to have be
considered earlier, with the best current result due to Iwaniec,
Kohnen and Sengupta~\cite{IKS07}. We will improve it, and obtain some
first results concerning the other two problems.  We will also suggest
further questions that may be of interest.

Before stating our main theorems, here are the basic notation about
modular forms (see, e.g.,~\cite[Ch. 14]{IK04} for a survey of these
facts). We denote by ${\rm H}_k^*(N)$ the finite set of all primitive
forms of weight $k$ for $\Gamma_0(N)$, where $k\ge 2$ is an even
integer and $N\geq 1$ is an integer. The restriction to trivial
Nebentypus ensures that all Fourier coefficients are real, and for any
$f\in {\rm H}_k^*(N)$, we denote
\begin{equation*}\label{SFf}
f(z) 
= \sum_{n=1}^\infty \lambda_f(n) n^{(k-1)/2} e(nz),\quad\quad
e(z)=e^{2i\pi z},\quad
\qquad(\Im m z>0),
\end{equation*}
its Fourier expansion at infinity. Since $f$ is primitive, the
$\lambda_f(n)$ are the normalized eigenvalues of the Hecke operators
$T_n$, and satisfy the well-known Hecke relations
\begin{equation}\label{Hecke}
\lambda_f(m) \lambda_f(n)
= \sum_{\substack{d\mid (m, n)\\ (d, N)=1}} 
\lambda_f\bigg(\frac{mn}{d^2}\bigg),
\end{equation}
for all integers $m\ge 1$ and $n\ge 1$.  In particular, $\lambda_f$ is
a multiplicative function of $n$ (so $\lambda_f(1)=1$) and moreover
the following important special case
\begin{equation}\label{Heckep}
\lambda_f(p)^2
= 1 + \lambda_f(p^{2})
\end{equation}
holds for all primes $p\nmid N$.
\par
Furthermore, it is also known that $\lambda_f(n)$ satisfies the deep
inequality
\begin{equation}\label{Deligne}
|\lambda_f(n)| \le \tau(n)
\end{equation}
for all $n\ge 1$, where $\tau(n)$ is the divisor function (this is the
Ramanujan-Petersson conjecture, proved by Deligne). In particular, we
have $\lambda_f(p)\in [-2,2]$ for $p\nmid N$, and hence there exists a
unique angle $\theta_f(p)\in [0,\pi]$ such that
\begin{equation}\label{eq-theta-f}
\lambda_f(p)=2\cos\theta_f(p).
\end{equation}
\par
Our other notation is standard in analytic number theory: for
instance, $\pi(x)$ denotes the number of primes $\leq x$ and
$P^+(n)$ (resp. $P^-(n)$) denotes the largest (resp. smallest) prime
factor of $n$, with the convention $P^+(1)=1$
(resp. $P^-(1)=\infty$). 
\par
We now describe our results.

\subsection{The first negative Hecke eigenvalue}

For $f\in \Hf_k^*(N)$, $k\geq 2$ and $N\geq 1$, it is well-known that
the coefficients $\lambda_f(n)$ change sign infinitely often. We
denote by $n_f$ the smallest integer $n\geq 1$ such that $(n,N)=1$ and
\begin{equation}\label{defnf}
\lambda_f(n)<0.
\end{equation}
\par
The analogue (or one analogue) of the least-quadratic non-residue
problem is to estimate $n_f$ in terms of the analytic conductor
$Q:=k^2N$.  Iwaniec, Kohnen and Sengupta \cite{IKS07} have shown
recently that
$$
n_f\ll Q^{29/60}=(k^2N)^{29/60}
$$
(here, standard methods lead to $n_f\ll_{\eps} Q^{1/2+\eps}$, so the
significance is that the exponent is $<1/2$).
\par
Our first result is a sharpening of this estimate:

\begin{theorem}\label{FirstNegative}
  Let $k\ge 2$ be an even integer and $N\ge 1$.  Then for all $f\in
  {\rm H}_k^*(N)$, we have
\begin{equation}\label{UB}
n_f\ll Q^{9/20}=(k^2N)^{9/20},
\end{equation}
where the implied constant is absolute.
\end{theorem}

This bound is not the best that can be achieved using our method, and
we will comment on this after its proof (in particular, an interesting
function $\beta$ occurs when trying to push the idea to its limit).
\par
One automatic improvement of the exponent arises from any subconvexity
bound for the relevant $L$-functions, as already observed by Iwaniec,
Kohnen and Sengupta. We do not need such deep results to prove
Theorem~\ref{FirstNegative}, but we will state below the precise
relation.
\par
We do not know if the estimate of Theorem~\ref{FirstNegative} holds
for the first negative Hecke eigenvalue at a prime argument.

\subsection{Statistic study of the first sign-change}

The upper bound (\ref{UB}) is probably far from optimal.  Indeed, one
can show that under the Grand Riemann Hypothesis we have
$$
n_f\ll (\log(kN))^2
$$
where the implied constant is absolute. Our next result confirms this
unconditionally for almost all $f$. It closely parallels the case of
Dirichlet characters (see \cite{LW08a}). Precisely, we first recall
that
$$
|\Hf_k^*(N)|\asymp k\varphi(N),
$$
where $\varphi(N)$ is the Euler function, as $k$, $N\ra +\infty$, and
we prove:

\begin{theorem}\label{FirstNegativeAlmost}
  Let $\nu\ge 1$ be a fixed integer and ${\mathscr P}$ be a set of
  prime numbers of positive density in the following sense:
$$
\sum_{\substack{z<p\le 2z\\ p\in {\mathscr P}}} \frac{1}{p}
\ge \frac{\delta}{\log z} \qquad(z\ge z_0)
$$ 
for some constants $\delta>0$ and $z_0>0$.  Let $\{\eps_p\}_{p\in
  {\mathscr P}}$ be a sequence of real numbers such that $|\eps_p|=1$
for all $p$.  Let $k\ge 2$ be an even integer and $N\ge 1$ be
squarefree.  Then there are two positive constants $C$ and $c$ such
that the number of primitive cusp forms $f\in {\rm H}_k^*(N)$
satisfying
$$
\eps_p\lambda_f(p^\nu)>0
\quad
\hbox{for}
\quad
p\in {\mathscr P},
\quad
p \nmid N
\quad\hbox{and}\quad 
C\log(kN)<p\le 2C\log(kN)
$$
is bounded by
$$
\ll_{\nu, {\mathscr P}} kN 
\expll{-c}{kN}.
$$
\par
Here $C, c$ and the implied constant depend on $\nu$ and ${\mathscr
  P}$ only.
\end{theorem}

Taking ${\mathscr P}$ the set of all primes, $\eps_p=1$ and $\nu=1$ in
Theorem \ref{FirstNegativeAlmost}, we immediately get:\footnote{\ The
  cases where $\nu\geq 2$ can be interpreted as similar statements for
  the $\nu$-th symmetric powers.}

\begin{corollary}\label{cor-statistic}
  Let $k\ge 2$ be an even integer and $N\ge 1$ be squarefree.  There
  is an absolute positive constant $c$ such that we have
$$
n_f\ll \log(kN),
$$
for all $f\in \Hf_k^*(N)$, except for $f$ in an exceptional set with
$$
\ll kN\expll{-c}{kN}
$$
elements, where the implied constants are absolute.
\end{corollary}

It is very natural to ask whether this result is optimal (as the
analogue is known to be for real Dirichlet characters). In this
direction, we can prove the following:

\begin{theorem}\label{th-lower-bound}
  Let $N$ be a squarefree number and $k\geq 2$ an even integer, and
  let $(\eps_p)$ be a sequence of signs indexed by prime numbers. For
  any $\eps>0$, $\eps<1/2$, there exists $c>0$ such that
$$
\frac{1}{|\Hf_k^*(N)|}
|\{f\in \Hf_k^*(N)\,\mid\, \lambda_f(p)\text{ has sign } \eps_p\text{ for }
p\leq z,\ p\nmid N\}|\geq \Bigl(\frac{1}{2}-\eps\Bigr)^{\pi(z)}
$$
for $z=c\sqrt{(\log kN)(\log\log kN)}$, provided $kN$ is large enough.
\end{theorem}

One may expect that the same result would be true for $z\leq c\log kN$
(note that
$$
\Bigl(\frac{1}{2}\Bigr)^{\pi(c\log kN)}\geq
\expll{-c_1}{kN}
$$
so this result would be quite close to the statistic upper-bound of
Corollary~\ref{cor-statistic}, and would essentially be best possible,
confirming that the signs of $\lambda_f(p)$ behave almost like
independent (and unbiased) random variables in that range of $p$).
\par
Theorems~\ref{FirstNegativeAlmost} and~\ref{th-lower-bound} will be
proved in Section~\ref{sec-statistic}, using the method in
\cite{LW08b} and quantitative equidistribution statements for Hecke
eigenvalues, respectively.

\subsection{Recognition of modular forms by signs of Hecke eigenvalues}

Here we consider whether it is true that a primitive form $f$ is
determined uniquely by the sequence of signs of its Fourier
coefficients $\lambda_f(p)$, where we recall that we interpret the
sign of $0$ in a relaxed way, so that $0$ has the same sign as both
positive and negative numbers. 
\par
The answer to this question is, indeed, yes, and in fact (in the
non-CM case) an analogue of the strong multiplicity one theorem holds:
not only can we exclude finitely many primes, or a set of primes of
density zero, but even a set of sufficiently small positive
density. Here, the density we use is the analytic density defined as
follows: a set $E$ of primes has density $\kappa>0$ if and only if
\begin{equation}\label{eq-density}
\sum_{p\in E}{\frac{1}{p^{\sigma}}}\sim \kappa
\sum_{p}{\frac{1}{p^{\sigma}}}\sim -\kappa\log(\sigma-1)
\qquad
(\sigma\to 1+).
\end{equation}
\par
We will prove:

\begin{theorem}\label{th-1}
  Let $k_1, k_2\geq 2$ be even integers, let $N_1, N_2\geq 1$ be
  integers and $f_1\in {\rm H}_{k_1}^*(N_1)$, $f_2\in {\rm
    H}_{k_2}^*(N_2)$. 
\par
\emph{(1)} If the signs of $\lambda_{f_1}(p)$ and $\lambda_{f_2}(p)$
are the same for all $p$ except those in a set of analytic density
$0$, then $f_1=f_2$.
\par
\emph{(2)} Assume that neither of $f_1$ and $f_2$ is of CM type, for
instance assume that $N_1$ and $N_2$ are squarefree. Then, if
$\lambda_{f_1}(p)$ and $\lambda_{f_2}(p)$ have same sign for every
prime $p$, except those in a set $E$ of analytic density $\kappa$,
with $\kappa\leq 1/32$,  it follows that $f_1=f_2$.
\end{theorem}

Recall that a form $f\in \Hf^*_k(N)$ is of CM type if there exists a
non-trivial primitive real Dirichlet character $\chi$ such that
$\lambda_f(p)=\chi(p)\lambda_f(p)$ for all but finitely many primes
$p$. In that case, $\lambda_f(p)=0$ for all $p$ such that
$\chi(p)=-1$, and hence its signs coincide (in our relaxed sense) with
those of any other modular form for a set of primes of density at
least $1/2$. 
\par
Of course, Theorem \ref{th-1} is also valid for the natural density,
since the existence of the latter implies that of the analytic
density, and that they are equal. As a corollary, we get of course:

\begin{corollary}
  For any sequence of signs $(\eps_p)$ indexed by primes, there is at
  most pair $(k,N)$ and one $f\in\Hf_k^*(N)$ such that $\lambda_f(p)$
  has sign $\eps_p$ for all primes.
\end{corollary}

Theorem~\ref{th-1} is proved in Section~\ref{sec-recognition}. The
argument is short and simple, but it depends crucially on a very deep
result: Ramakrishnan's proof~\cite{ramakrishnan} the Rankin-Selberg
convolution $L$-function is the $L$-function of some modular form on
$GL(4)$.

\subsection{Motivation, further remarks and problems}

The main remark is that, underlying most of the problems we consider
is the Sato-Tate conjecture, which we recall (see Mazur's
survey~\cite{mazur}): provided $f$ is not of CM type (for instance, if
$N$ is squarefree), one should have
$$
\lim_{x\ra +\infty}{\frac{1}{\pi(x)} |\{p\leq x\,\mid\, \theta_f(p)\in
  [\alpha,\beta]\}|}= \int_{\alpha}^{\beta}{\d\mu_{ST}}
$$
for any $\alpha<\beta$, where $\mu_{ST}$ is the Sato-Tate measure
$$
\mu_{ST}=\frac{2}{\pi}\sin^2\theta \d \theta,
$$
on $[0,\pi]$. Since $\mu_{ST}([0,\pi/2])=\mu_{ST}([\pi/2,\pi])$, this
indicates in particular that the signs of $\lambda_f(p)$ should be
equitably shared between $+1$ and $-1$. This suggests and motivates
many of our results and techniques of proof.
\par
There is much ongoing progress on the Sato-Tate conjecture; for $f\in
\Hf_k^*(N)$, non-CM, a proof of the conjecture has been announced by
Barnet-Lamb, Geraghty, Harris and
Taylor~\cite[Th. B]{satotate}. However, knowing its truth does not
immediately simplify our arguments. Indeed, it would be most
immediately relevant for parts of Section~\ref{sec-recognition}, but
Theorem~\ref{th-1} is really concerned with the \emph{pair-Sato-Tate
  conjecture}, or property, which would be the statement, for a pair
$(f_1,f_2)$, that for any $a_1<b_1$, $a_2<b_2$ the set of primes
$$
\big\{p\,\mid\, \lambda_{f_1}(p)\in [a_1,b_1]\text{ and }
\lambda_{f_2}(p)\in [a_2,b_2]\big\}
$$
has density equal to $\mu_{ST}([a_1,b_1])\mu_{ST}([a_2,b_2])$ (in
other words, the Fourier coefficients at primes are
\emph{independently} Sato-Tate distributed). This is expected to hold
for any pair of non-CM modular forms, such that neither is a quadratic
twist of the other (and in the case of elliptic curves,
Mazur~\cite[Footnote 12]{mazur} mentions that there is ongoing
progress by Harris on this problem). If this holds, it will follow
that the density of coincidences of signs is always $\leq \dm$,
which is the probability under independent Sato-Tate measures that two
``samples'' are of the same sign.
\par
Here is a natural question which is suggested by Theorem~\ref{th-1}:
estimate the size, as a function of the weight and conductor, of the
smallest integer $n_{f_1,f_2}$ for which the sign of
$\lambda_{f_1}(n)$ and $\lambda_{f_2}(n)$ are different. If we enlarge
slightly our setting to allow $f_2$ to be an Eisenstein series such as
the Eisenstein series of weight $4$:
$$
E_4(z)=1+240\sum_{n\geq 1}{\Bigl(\sum_{d\mid n}{d^3}\Bigr)e(nz)},
$$
where all Hecke eigenvalues are positive, then the question becomes
(once more) that of finding the first negative Hecke eigenvalue for
$f_1$, i.e., the problem considered in
Theorem~\ref{FirstNegative}. Hence, we know that
$$
n_{f_1,E_4}\ll (k_1^2N_1)^{10/21},
$$
where the implied constant is absolute, but it would be interesting to
obtain a more general version, in particular a uniform one with
respect to both $f_1$ and $f_2$.
\par
At least our statistic result (Theorem~\ref{FirstNegativeAlmost})
generalizes immediately if one of the forms is fixed: taking
${\mathscr P}$ to be the set of all primes, $\nu=1$ and
$\eps_p={\rm sign}\,\lambda_{f_2}(p)$ if $\lambda_{f_2}(p)\not=0$, and
$1$ otherwise, we get 
immediately the following corollary: 

\begin{corollary}
\label{FirstSign}
  Let $k_1, k_2\geq 2$ be even integers and $N_1, N_2\geq 1$
  squarefree.  For any fixed $f_2\in {\rm H}_{k_2}^*(N_2)$, there is
  an absolute positive constant $c$ such that 
$$
n_{f_1, f_2}\ll_{f_2} \log(k_1N_1),
$$
for all $f\in \Hf_{k_1}^*(N_1)$ except for those in an exceptional set
with
$$
\ll k_1N_1\expll{-c}{k_1N_1}
$$
elements, where the implied constants depend only on $f_2$.
\end{corollary}


\section{Proof of Theorem \ref{FirstNegative}}
\label{sec-firstneg}

Although some ideas are related to those of Iwaniec, Kohnen and
Sengupta, the proof of Theorem~\ref{FirstNegative} is somewhat easier.
\par
Thus let $f\in H_k^*(N)$, and let $y>0$ be such that $\lambda_f(n)\geq
0$ for $n\le y$ and $(n, N)=1$.  The idea to estimate $y$ is to
compare upper and lower bounds for the sum
$$
S(f,x)=\mathop{{\sum_{\stacksum{n\le x}{(n,N)=1}}}^{\hskip -2,5mm \flat}} \hskip 1,5mm {\lambda_f(n)},
$$
where $\sum^\flat$ restricts a sum to squarefree integers.  An upper
bound is easily achieved: using the convexity bound for Hecke
$L$-functions and the Perron formula, we obtain
\begin{equation}\label{eq-sf-upper}
S(f,x)\ll_\varepsilon (k^2N)^{1/4+\varepsilon} x^{1/2+\varepsilon}
\quad(x\ge 1).
\end{equation}
\par
This estimate is independent of any information on $y$. We note that,
more generally, if we have
$$
L(f,1/2+it)\ll (k^2N(1+|t|)^2)^{\eta},
$$
for $t\in\R$, where $\eta>0$, then we get~(\ref{eq-sf-upper}) for
$x\geq Q^{2\eta+\eps}$ (the recent work of Michel and
Venkatesh~\cite{mv} provides such a uniform result for some -- very
small -- $\eta<1/4$, and the convexity bound states that any
$\eta>1/4$ is suitable).
\par
We now proceed to establish a lower bound for $S(f,x)$ by using the
assumption. For primes $p\le y$ with $p\nmid N$, we thus have
$\lambda_f(p)\geq 0$; furthermore, if $p\le \sqrt{y}$ and $p\nmid N$,
we have the better bound $\lambda_f(p)\geq 1$ because
$$
\lambda_f(p)^2=1+\lambda_f(p^2)\geq 1
$$ 
by the Hecke relation.
\par
We now introduce an auxiliary multiplicative function $h=h_y$ defined
by
$$
h_y(p)=\begin{cases}
-2    & \text{if $\,p>y$ and $p\nmid N$}
\\
0     & \text{if $\,\sqrt{y}<p\le y$ or $\,p\mid N$}
\\
1     & \text{if $\,p\le \sqrt{y}$ and $p\nmid N$}

\end{cases}
$$
and $h_y(p^\nu)=0$ for $\nu\ge 2$.  
\par
We shall show in a moment the following lemma:

\begin{lemma}\label{lm-h}
  For any $y>0$, define $h_y$ as above. Then, for any $\varepsilon>0$,
  we have
\begin{equation}\label{hyn0}
\sum_{n\le y^u } h_y(n)
= \zeta_M(2)^{-1} \frac{\varphi(N)}{N}y^u \big(\rho(2u)-2\log u\big)
\bigg\{1+O\bigg(\frac{(\log_2y)^2}{\log y}\bigg)\bigg\}
\end{equation} 
uniformly for
\begin{equation}\label{ConditonuyN}
  1\le u\le \tfrac{3}{2}
  \qquad\text{and}\qquad
  y\ge N^{1/3},
\end{equation}
where 
$$
\zeta_N(2)=\prod_{p\nmid N}(1-p^{-2})^{-1}
$$ 
and $\rho(u)$ is the Dickman function, defined as the unique continous
solution of the difference-differential equation
$$
u\rho'(u)+\rho(u-1)=0
\quad(u>1),
\qquad
\rho(u)=1
\quad(0<u\le 1).
$$
In particular $\rho(2u)-2\log u>0$ for all $u < \kappa$ where $\kappa$ is
the solution to $\rho(2\kappa)=2\log \kappa$. 
We have $\kappa>\tfrac{10}{9}$.
\end{lemma}
The point of introducing this auxiliary function is that, with
notation as before, we have the lower bound
\begin{equation}\label{pointdepart}
S(f,y^u)\ge \sum_{n\le y^u} h_y(n)
\end{equation}
for all $u< \kappa$, provided $y$ is large enough, e.g., $y\geq
N^{1/3}$ with $N$ large enough, which we can obviously assume in
proving Theorem~\ref{th-1}.
\par
Indeed, let $g_y$ be the multiplicative function defined by the
Dirichlet convolution identity
$$
\lambda_f=g_y*h_y,
$$
then $g_y(n)\geq 0$ for all squarefree integers $n\geq 1$ such that
$(n, N)=1$, since
$$
g_y(p)=\lambda_f(p)-h_y(p)\geq 0
$$
for all $p\nmid N$ (the case $p\geq y$ following from Deligne's
inequality). We have trivially
$$
\sum_{n\leq z}{h_y(n)}\geq 0
$$
if $z\leq y$ since each term is non-negative in that range, and
additionally
$$
\sum_{n\leq y^u}{h_y(n)}\geq 0
$$
for $u<\kappa$ if $y$ is large enough, by Lemma~\ref{lm-h}.
\par
Hence
$$
S(f,y^u)=\sumb_{\stacksum{n\leq y^u}{(n,N)=1}} {\lambda_f(n)}
=\sumb_{d\leq y^u}{g_y(d)\sumb_{\ell\leq y^u/d}{h_y(\ell)}} \\
\geq
\sumb_{\ell\leq y^u}{h_y(\ell)}
$$
since every term in the sum over $d$ is non-negative and $g_y(1)=1$.
\par
We now deduce from Lemma~\ref{lm-h}
$$
S(f,y^u) \geq \sum_{n\le y^u} h_y(n) \gg \frac{y^u}{\log\log N},
$$
for $u<\kappa$. Then, a comparison with~(\ref{eq-sf-upper}) gives the
estimate
$$ 
y\le (k^2 N)^{1/(2\kappa) + o(1)}. 
$$
\par
Quoting the lower bound for $\kappa$ from Lemma~\ref{lm-h}, we are done.

\begin{proof}[Proof of Lemma~\ref{lm-h}]
According to the definition of $h_y$, we have
\begin{equation}\label{hyn1} 
\sum_{n\le y^u} h_y(n) 
= \sumb_{\substack{n\le y^u\\ P(n)\le \sqrt{y}}} 1 
- 2 \sum_{\substack{y < p\le y^u\\ p\nmid N}} \sumb_{m\le y^u/p} 1
\end{equation}
for all $u$ and $y$ satisfying \eqref{ConditonuyN}, where we use the
convention here that $\sumb$ also restricts the sum to $(n,N)=1$ to
simplify notation.
\par
The second term contributes
$$
-\frac{2}{\zeta_N(2)}\frac{\varphi(N)}{N}y^u(\log u)\bigg\{1+O\biggl(
\frac{(\log_2 y)^2}{\log y}\biggr)\bigg\}
$$
by standard estimates (the leading constant arises of course because
$$
\sumb_{\stacksum{n\le
    y^u}{(n,N)=1}}{1}\sim \frac{1}{\zeta_N(2)}\frac{\varphi(N)}{N}y^u
$$
uniformly in our range).
\par
For the first term, if it were not for this condition and the
requirement that $n$ be squarefree, the lemma would then follow
immediately from the well-known property
$$
\sum_{\substack{n\le y^u\\ P(n)\le \sqrt{y}}} 1= \rho(2u)y^u
\bigg\{1+O\bigg(\frac{1}{\log y}\bigg)\bigg\},
$$
of the Dickman function. However, because our uniformity requirements
are quite modest, it is fairly simple to deduce the stated inequality
from this result using Möbius inversion to detect the coprimality and
squarefree condition (for very general bounds of this type,
see~\cite[Th. 2.1]{tw}, though our requirements are much weaker).
\par
A numerical computation using \textsc{Maple} leads to
$\kappa>\tfrac{10}{9}$.
\end{proof} 

\begin{remark}
  To estimate $\kappa$, one can also use the lower bound
$$
\rho(2u)
= 1- \log (2u) + \int_2^{2u} \frac{\log(t-1)}{t} \d t
\ge 1- \log (2u)
\qquad(1\le u\le \tfrac{3}{2}),
$$ 
which shows that $\rho(2u)-2\log u\ge 1- \log 2 -3 \log u$, and
$$ 
\kappa \ge (e/2)^{1/3}>\tfrac{11}{10},
$$ 
which leads to the exponent $\frac{5}{11}=0.4545\ldots$ in
Theorem~\ref{FirstNegative}.
\end{remark}

\begin{remark}
  This result is not the limit of the method employed. Precisely, in
  addition to obtained $\lambda_f(p)\geq 1$ for $p<\sqrt{y}$,
  $(p,N)=1$, we can exploit higher powers: write
  $\lambda_f(p)=2\cos\theta_f(p)$ with $\theta_f(p)\in [0,\pi$]. Then
  if $m\geq 1$ is an integer, we have for $1\leq j\leq m$ and
  $y^{1/(m+1)}\leq p<y^{1/m}$ that
$$
0\leq \lambda_f(p^j)=\frac{\sin ((j+1)\theta_f(p))}{\sin\theta_f(p)}
$$
(if $p\nmid N$). This implies that $\theta_f(p)\leq \pi/(m+1)$, and
hence 
$$
\lambda_f(p)\geq 2\cos\frac{\pi}{m+1},
$$
for $p\nmid N$ with $y^{1/(m+1)}\leq p<y^{1/m}$. This can be exploited
by boundind $S(f,y)$ from below using a new auxiliary function $h$
supported on squarefree numbers coprime to $N$ with
$$
h(p)=\alpha\Bigl(\frac{\log p}{\log y}\Bigr)
$$
where $\alpha(u)=-2$ for $u\geq 1$, $\alpha(0)=2$ and
$\alpha(u)=2\cos(\pi/(m+1))$ if $1/(m+1)\leq u<1/m$. One can show an
asymptotic of the type
$$
\sum_{n\leq y^u}{h(n)}\sim C\beta(u)y^u(\log y^u)
$$
(assuming $N=1$ for simplicity) for some constant $C>0$, where the
function $\beta$ can be described by the following inclusion-exclusion
formula:
$$
u\beta(u)=u+\sum_{j\geq 1}{\frac{(-1)^j}{j!}I_j(u)},
$$
with
$$
I_j(u)=\int_{\Delta_j}{(u-t_1-\cdots-t_j)\prod_{i=1}^j{(2-\alpha(t_j))}
\frac{dt_1\cdots dt_j}{t_1\cdots t_j}},
$$
integration ranging over the set
$$
\Delta_j=\{(t_1,\ldots,t_j)\in [0,+\infty[^j\,\mid\,
t_1+\cdots+t_j\leq u\}.
$$
\par
This function is also a solution of the integral equation
$$
u^2\beta(u)=\int_0^u{t\beta(t)\alpha(u-t)dt}
$$
(see~\cite{gs} for related investigations of a class of integral
equations of this type). 
\par
To improve Theorem~\ref{FirstNegative}, one needs to find (a close
approximation to) the first positive zero of $\beta$. We have not
found a nice way to compute $\beta$ numerically, but this would be
quite an interesting problem, and its solution is likely to lead to
significant improvements in the result. We hope to come back to this
in the future.
\end{remark}

\section{Statistical results}
\label{sec-statistic}

Our goal is now to prove Theorems~\ref{FirstNegativeAlmost}
and~\ref{th-lower-bound}. For the first, the main tool is the
following type of large sieve inequality.

\begin{lemma}[\cite{LW08b}, Theorem 1]\label{LargeSieve}
Let $\nu\ge 1$ be a fixed integer
and let $\{b_p\}_p$ be a sequence of real numbers
indexed by prime numbers such that $|b_p|\le B$ 
for some constant $B$ and for all primes $p$.
Then we have
$$
\sum_{f\in {\rm H}_k^*(N)} 
\bigg|\sum_{\substack{P<p\le Q\\ p\,\nmid N}} 
b_p \frac{\lambda_f(p^\nu)}{p}\bigg|^{2j}
\ll_\nu k\varphi(N) 
\bigg(\frac{96B^2(\nu+1)^2j}{P\log P}\bigg)^j 
+ (kN)^{10/11}\bigg(\frac{10BQ^{\nu/10}}{\log P}\bigg)^{2j}
$$
uniformly for 
$$B>0,
\qquad
j\ge 1,
\qquad
2\mid k, 
\qquad
2\le P<Q\le 2P,
\qquad
N\ge 1 \quad (\hbox{squarefree}).
$$
The implied constant depends on $\nu$ only.
\end{lemma}

\begin{proof}[Proof of Theorem~\ref{FirstNegativeAlmost}]
The basic idea is that for all forms $f$ with coefficients
$\lambda_f(p^{\nu})$ of the same sign $\eps_p$, the sums
$$
\sum_{\substack{P<p\le 2P\\ p\in {\mathscr
      P}}}{\frac{\eps_p\lambda_f(p^{\nu})}{p}} 
$$
exhibit no cancellation due to variation of signs. The large sieve
implies this is very unlikely to happen, except if the
$\lambda_f(p^{\nu})$ are very small in absolute value. The Hecke
relations are used to control this other possibility by relating it to
$\lambda_f(p^{2\nu})$ being large which can not happen too often
either.\footnote{\ Variants of this well-known trick have been used in
  a number of other contexts, as in~\cite{DK00}, but note that the
  large sieve inequality proved there would not work for this problem,
  due to the lack of multiplicative stability of the sign conditions
  (it would also be much less efficient).}
\par
For the details, we first denote
$$
{\mathscr P}_N:=\{p\in {\mathscr P}\,\mid\,p\nmid N\},
$$
and define
\begin{align*}
{\mathscr E}_k^*(N, P; {\mathscr P}) 
& := \big\{ f\in {\rm H}_k^*(N)\,\mid\,
\eps_p \lambda_f(p^\nu)>0 \;\; \hbox{for} \;\;
p\in {\mathscr P}_N\cap (P, 2P]\big\},
\\\noalign{\vskip 1mm}
{\mathscr E}_k^{\nu'}(N, P; {\mathscr P})
& := \bigg\{ f\in {\rm H}_k^*(N)\,\mid\,
\bigg|\sum_{\substack{P<p\le 2P\\ p\in {\mathscr P}_N}} 
\frac{\lambda_f(p^{2\nu'})}{p}\bigg|
\ge \frac{\delta}{2\nu\log P}\bigg\}
\quad(1\le \nu'\le \nu).
\end{align*}
\par
To prove Theorem~\ref{FirstNegativeAlmost}, clearly we only need to
show that there are two positive constants $C=C(\nu, {\mathscr P})$
and $c=c(\nu, {\mathscr P})$ such that
\begin{equation}\label{pthm1a}
|{\mathscr E}_k^*(N, P; {\mathscr P})|
\ll_{\nu, {\mathscr P}} kN \expll{-c}{kN}
\end{equation}
uniformly for
$$2\mid k,
\qquad
N\quad({\rm squarefree}),
\qquad
kN\ge X_0,
\qquad
C\log(kN)\le P\le (\log(kN))^{10}
$$
for some sufficiently large number $X_0=X_0(\nu, {\mathscr P})$.
\par
The definition of ${\mathscr E}_k^*(N, P; {\mathscr P})$ and Deligne's
inequality allow us to write
\begin{align*}
\sum_{f\in {\mathscr E}_k^*(N, P; {\mathscr P})} 
\bigg|\sum_{\substack{P<p\le 2P\\ p\in {\mathscr P}_N}} 
\frac{\lambda_f(p^\nu)^2}{p}\bigg|^{2j}
& \le \sum_{f\in {\mathscr E}_k^*(N, P; {\mathscr P})} 
\bigg|\sum_{\substack{P<p\le 2P\\ p\in {\mathscr P}_N}} 
(\nu+1) \eps_p \frac{\lambda_f(p^\nu)}{p}\bigg|^{2j}
\\
& \le \sum_{f\in {\rm H}_k^*(N)} 
\bigg|\sum_{\substack{P<p\le 2P\\ p\in {\mathscr P}_N}} 
(\nu+1)\eps_p \frac{\lambda_f(p^\nu)}{p}\bigg|^{2j}.
\end{align*}
Choosing
$$
b_p = \begin{cases}
(\nu+1)\eps_p  & \text{if $p\in {\mathscr P}$},
\\
0                   & \text{otherwise}
\end{cases}
$$
in Lemma \ref{LargeSieve}, we find that
\begin{equation}\label{pthm1c}
\begin{aligned}
\sum_{f\in {\mathscr E}_k^*(N, P; {\mathscr P})} 
\bigg|\sum_{\substack{P<p\le 2P\\ p\in {\mathscr P}_N}} 
\frac{\lambda_f(p^\nu)^2}{p}\bigg|^{2j}
& \le \sum_{f\in {\rm H}_k^*(N)} 
\bigg|\sum_{\substack{P<p\le 2P\\ p\,\nmid N}} 
b_p \frac{\lambda_f(p^\nu)}{p}\bigg|^{2j}
\\
& \ll kN \bigg(\frac{96(\nu+1)^4j}{P\log P}\bigg)^j 
+ (kN)^{10/11} P^{\nu j/2}.
\end{aligned}
\end{equation} 
\par
In view of the Hecke relation (\ref{Hecke}), the left-hand side of
(\ref{pthm1c}) is
\begin{align*}
& \ge \sum_{f\in {\mathscr E}_k^*(N, P; {\mathscr P})\setminus 
(\cup_{\nu'=1}^\nu{\mathscr E}_k^{\nu'}(N, P; {\mathscr P}))} 
\bigg(
\sum_{\substack{P<p\le 2P\\ p\in {\mathscr P}_N}} \frac{1}{p}
- \sum_{1\le \nu'\le \nu} \bigg|
\sum_{\substack{P<p\le 2P\\ p\in {\mathscr P}_N}} 
\frac{\lambda_f(p^{2\nu'})}{p}
\bigg|
\bigg)^{2j}
\\
& \ge \sum_{f\in {\mathscr E}_k^*(N, P; {\mathscr P})\setminus 
(\cup_{\nu'=1}^\nu{\mathscr E}_k^{\nu'}(N, P; {\mathscr P}))} 
\bigg(
\sum_{\substack{P<p\le 2P\\ p\in {\mathscr P}_N}} \frac{1}{p}
- \frac{\delta}{2\log P}\bigg)^{2j}.
\end{align*} 
\par
Let $\omega(n)$ be the number of distinct prime factors of $n$.
Using the hypothesis on ${\mathscr P}$ and the classical inequality
$$\omega(n)\le \{1+o(1)\}\frac{\log n}{\log\log n},$$ 
we infer that
\begin{align*}
\sum_{\substack{P<p\le 2P\\ p\in {\mathscr P}_N}} \frac{1}{p}
- \frac{\delta}{2\log P}
& \ge 
\sum_{\substack{P<p\le 2P\\ p\in {\mathscr P}}} \frac{1}{p}
- \sum_{\substack{P<p\le 2P\\ p\,\mid N}} \frac{1}{p}
- \frac{\delta}{2\log P}
\\&
 \ge \frac{\delta}{2\log P} - \frac{\omega(N)}{P}
\ge \frac{\delta/2-2/C}{\log P}
 \ge \frac{\delta}{6\log P},
\end{align*}
provided $C\ge 6/\delta$. Combining this with (\ref{pthm1c}), we infer
that
$$
|{\mathscr E}_k^*(N, P; {\mathscr P})\setminus
(\cup_{\nu'=1}^\nu{\mathscr E}_k^{\nu'}(N, P; {\mathscr P}))| \ll kN
\bigg(\frac{3456(\nu+1)^4j\log P}{\delta^2P}\bigg)^j + (kN)^{10/11}
P^{j}.
$$
\par
Now we bound the size of the sets ${\mathscr E}_k^{\nu'}(N, P;
{\mathscr P})$ to finish the proof. Taking
$$
B=1,
\qquad 
\nu=2\nu',
\qquad 
Q=2P
\qquad{\rm and}\qquad
b_p = \begin{cases}
1 & \text{if $p\in {\mathscr P}$}
\\
0 & \text{otherwise}
\end{cases}
$$ 
in Lemma \ref{LargeSieve}, we get
\begin{align*}
\bigg(\frac{\delta}{2\log P}\bigg)^{2j} 
|{\mathscr E}_k^{\nu'}(N, P; {\mathscr P})|
& \le \sum_{f\in {\rm H}_k^*(N)} 
\bigg|\sum_{\substack{P<p\le 2P\\ p\,\nmid N}} 
b_p \frac{\lambda_f(p^{2\nu'})}{p}\bigg|^{2j}
\\
& \ll kN \bigg(\frac{96(2\nu'+1)^2j}{P\log P}\bigg)^j 
+  (kN)^{10/11}\bigg(\frac{10(2P)^{\nu'/5}}{\log P}\bigg)^{2j}.
\end{align*}
Hence, 
\begin{equation}\label{pthm1b}
|{\mathscr E}_k^{\nu'}(N, P; {\mathscr P})|
\ll kN \bigg(\frac{3456\nu^4j\log P}{\delta^2P}\bigg)^j 
+ (kN)^{10/11} P^{\nu j}
\quad(1\le \nu'\le \nu)
\end{equation} 
provided $P\ge 2(20\nu/\delta)^{10/(3\nu)}$. 
\par
Combining this with (\ref{pthm1b}), we finally obtain
\begin{equation}\label{pthm1d}
|{\mathscr E}_k^*(N, P; {\mathscr P})|
\ll kN \bigg(\frac{3456(\nu+1)^4j\log P}{\delta^2P}\bigg)^j 
+ (kN)^{10/11}  P^{\nu j}
\end{equation}
uniformly for
$$2\mid k,
\qquad
N\quad({\rm squarefree}),
\qquad
C\log(kN)\le P\le (\log(kN))^{10},
\qquad
j\ge 1.$$
\par
Now, take
$$
j = \bigg[\delta^* \frac{\log(kN)}{\log P}\bigg]
$$
where $\delta^*=\delta^2/(10(\nu+1))^4$.  We can ensure $j>1$ once
$X_0$ is chosen to be suitably large.  A simple computation gives that
$$
\bigg(\frac{3456(\nu+1)^4 j\log P}{\delta^2 P}\bigg)^j \ll
\expll{-c}{kN}
$$ 
for some positive constant $c=c(\nu, {\mathscr P})$ and $P^{\nu j} \ll
(kN)^{1/1000}, $ provided $X_0$ is large enough.  Inserting them into
(\ref{pthm1d}), we get (\ref{pthm1a}) and complete the proof.
\end{proof}

We now come to the lower bound of Theorem~\ref{th-lower-bound}.  Our
basic tool here is an equidistribution theorem for Hecke eigenvalues
which is of some independent interest: it shows (quantitatively) that,
after suitable average over $\Hf_k^*(N)$, the Hecke eigenvalues
corresponding to the first primes are independently Sato-Tate
distributed (thus, it is related to the earlier work of
Sarnak~\cite{sarnak} for Maass forms and Serre~\cite{serre} and
Royer~\cite{royer} for holomorphic forms).
\par
First, we recall the definition~(\ref{eq-theta-f}) of the angle
$\theta_f(p)\in [0,\pi]$ associated to any $f\in \Hf_k^*(N)$ and prime
$p\nmid N$. We also recall
that the Chebychev \emph{functions} $X_n$, $n\geq 0$, defined by
\begin{equation}\label{eq-cheby-pol}
X_n(\theta)=\frac{\sin((n+1)\theta)}{\sin\theta}
\end{equation}
for $\theta\in [0,\pi]$, form an orthonormal basis of
$L^2([0,\pi],\mu_{ST})$. Hence, for any $\omega\geq 1$, the functions
of the type
$$
(\theta_1,\ldots,\theta_{\omega})\mapsto \prod_{1\leq j\leq \omega}
{X_{n_j}(\theta_j)}
$$
for $n_j\geq 0$, form an orthonormal basis of
$L^2([0,\pi]^{\omega},\mu_{ST}^{\otimes \omega})$.

\begin{proposition}\label{pr-equi}
  Let $N$ be a squarefree number, $k\geq 2$ an even integer, $s\geq 1$
  an integer and $z\geq 2$ a real number. For any prime $p\leq z$
  coprime with $N$, let
$$
Y_p(\theta)=\sum_{j=0}^s{\hat{y}_{p}(j)X_j(\theta)}
$$
be a ``polynomial'' of degree $\leq s$ expressed in the basis of
Chebychev functions on $[0,\pi]$. Then we have
$$
\sum_{f\in \Hf_k^*(N)}{\omega_f\prod_{\stacksum{p\leq z}{(p,N)=1}}{
Y_p(\theta_f(p))
}}=
\prod_{\stacksum{p\leq z}{(p,N)=1}}{\hat{y}_p(0)}+
O(C^{\pi(z)}D^{sz}(\tau(N)\log 2N)^2(Nk^{5/6})^{-1})
$$
where 
\begin{gather*}
\omega_f=\frac{\Gamma(k-1)}{(4\pi)^{k-1}\langle f,f\rangle}
\frac{N}{\varphi(N)},\quad
\langle f,f\rangle\text{ the Petersson norm of $f$}, 
\\
C=\max_{p,j}{|\hat{y}_p(j)|}, 
\end{gather*}
and $D\geq 1$ and the implied constant are absolute.
\end{proposition}

By linearity, clearly, we get an analogue result for
$$
\sum_{f\in \Hf_k^*(N)}{\omega_f
\varphi((\theta_p)_{p\leq z})},\quad\quad
\varphi=\sum_{j}{\varphi_j},
$$
where each $\varphi_j$ is a function which is a product of polynomials
as in the statement.

\begin{proof}
  Using the fact that for any $n_p\geq 0$, we have
\begin{equation}\label{eq-cheby}
\prod_{\stacksum{p\leq z}{(p,N)=1}}{X_{n_p}(\theta_f(p))}=
\lambda_f\Bigl(\prod_{\stacksum{p\leq z}{(p,N)=1}}{p^{n_p}}\Bigr), 
\end{equation}
(which is another form of the Hecke multiplicativity), we expand the
product and get
$$
\prod_{\stacksum{p\leq z}{(p,N)=1}}{Y_p(\theta_f(p))}
=\sum_{d\mid P_N(z)^s}{
\Bigl(\prod_{p\mid P_N(z)}{\hat{y}_p(v_p(d))}\Bigr)\lambda_f(d)
}
$$
where $v_p(d)$ is the $p$-adic valuation of an integer and $P_N(z)$ is
the product of the primes $p\leq z$, $p\nmid N$.
\par
We now sum over $f$ and appeal to the following Petersson formula for
primitive forms:
$$
\sum_{f\in \Hf_k^*(N)}{\omega_f\lambda_f(m)}=
\delta(m,1)+O(m^{1/4}\tau(N)^2(\log 2mN)^2(Nk^{5/6})^{-1}),
$$
for all $m\geq 1$ coprime with $N$ (this is a simplified version of
that in~\cite[Cor. 2.10]{ILS}; note our slightly different definition
of $\omega_f$, which explains the absence of $\varphi(N)/N$ on the
right-hand side); the result then follows easily from simple estimates
for the sum over $d$ of the remainder terms.
\end{proof}

We now deduce Theorem~\ref{th-lower-bound} from this, assuming
$\eps_p=1$ for all $p$ (handling the other choices of signs being
merely a matter of complicating the notation). 
\par
To simplify notation, we write $P=P_N(z)$ the product of primes $\leq
z$ coprime with $N$, and $\omega$ the number of such primes.
\par
First, if we wanted only to have $\lambda_f(p)\geq 0$ for a fixed
(finite) set of primes (i.e., for $z$ fixed), we would be immediately
done: Proposition~\ref{pr-equi} shows\footnote{\ For Maass forms, this
  is essentially one of the early results of Sarnak~\cite{sarnak}.}
that the $(\theta_f(p))_{p\mid P}$ become equidistributed as $kN\ra
+\infty$ with respect to the product Sato-Tate measure, if we weigh
modular forms with $\omega_f$, and hence
$$
\sum_{\stacksum{f\in\Hf_k^*(N)}{p\mid P\Rightarrow \lambda_f(p)\geq
    0}}{\omega_f}\ra
\mu_{ST}([0,\pi/2])^{\omega}=\Bigl(\frac{1}{2}\Bigr)^{\omega}
$$
which is of the desired type, except for the presence of the
weight.\footnote{\ Using the trace formula instead of the Petersson
  formula (as in~\cite{royer}), the unweighted analogue of
  Proposition~\ref{pr-equi} holds with a product of local Plancherel
  measures, but each still gives measure $1/2$ to the two signs.}
However, we want to have
$$
\lambda_f(p)\geq 0\text{ for } p\leq z,\ (p,N)=1,
$$
where $z$ grows with $kN$, and this involves quantitative lower bounds
for approximation in large dimension, which requires more care. We use
a result of Barton, Montgomery and Vaaler~\cite{BMV00} for this
purpose; although it is optimized for uniform distribution modulo $1$
instead of the Sato-Tate context, but it is not difficult to adapt it
here and this gives a quick and clean argument.\footnote{\ This result
  was also used recently by Y. Lamzouri~\cite[\S 7]{L08}, in a
  somewhat related context.}
\par
Precisely, we consider $[0,\pi]^{\omega}$, with the product Sato-Tate
measure, and we will write $\theta=(\theta_p)$ for the elements of
this set; we also consider $[0,1]^{\omega}$ and we write $x=(x_p)$ for
elements
there. 
\par 
For any positive odd integer $L$, we get from~\cite[Th. 7]{BMV00} two
explicit trigonometric polynomials\footnote{\ Meaning, standard
  trigonometric polynomials of the type
  $\sum_{\ell}{\alpha_{\ell}e(\ell\cdot x)}$.} on $[0,1]^{\omega}$,
denoted $A_L(x)$, $B_L(x)$, such that
$$
A_L(\theta/\pi)-B_L(\theta/\pi)\leq \prod_{\stacksum{p\leq
    z}{(p,N)=1}}{ \chi(\theta_p) }
$$
for all $\theta=(\theta_p)\in [0,\pi]^{\omega}$, where
$\chi(\theta_p)$ is the characteristic function of $[0,\pi/2]\subset
[0,\pi]$ (precisely, we consider the functions denoted $A(x)$, $B(x)$
in~\cite{BMV00}, with parameters $N=\omega$ and $u_n=0$, $v_n=1/2$ for
all $n\leq \omega$; since $(v_n-u_n)(L+1)=(L+1)/2$ is a positive
integer, we are in the situation $\Phi_{u,v}\in \mathcal{B}_N(L)$ of
loc. cit.).
\par
Thus we have the lower bound
\begin{equation}\label{eq-th-3-lower}
\sum_{\stacksum{f\in \Hf_k^*(N)}{\lambda_f(p)\geq 0\text{ for }p\mid
    P}}{\omega_f }
\geq \sum_{f\in \Hf_k^*(N)}{\omega_f \Bigl(
A_L(\theta_f/\pi)-B_L(\theta_f/\pi)
\Bigr)},
\end{equation}
where $\theta_f=(\theta_f(p))_{p}$.
\par
Moreover, as we will explain below, $A_L(\theta/\pi)$ is a product of
polynomials over each variable, and $B_L(\theta/\pi)$ is a sum of
$\omega$ such products, and we can now apply Proposition~\ref{pr-equi}
(and the remark following it) to the terms on the right-hand
side. More precisely, we claim that the following lemma holds:

\begin{lemma}\label{lm-polys}
With notation as above, we have:
\par
\emph{(1)} For any $\eps\in (0,1/2)$, there exists constants $L_0\geq
1$, and $c>0$, such that the contribution $\Delta$ of the constant
terms of the Chebychev expansions of $A_L(\theta/\pi)$ and
$B_L(\theta/\pi)$ satisfies
$$
\Delta \geq \Bigl(\frac{1}{2}-\eps\Bigr)^{\pi(z)},
$$
if $L$ is the smallest odd integer $\geq c\pi(z)$ and if $L\geq L_0$.
\par
\emph{(2)} All the coefficients in the expansion in terms of Chebychev
functions of the factors in $A_L(\theta/\pi)$ or in the terms of
$B_L(\theta/\pi)$ are bounded by $1$.
\par
\emph{(3)} The degrees, in terms of Chebychev functions, of the
factors of $A_L(\theta/\pi)$ and of the terms of $B_L(\theta/\pi)$,
are $\leq 2L$.
\end{lemma}

Using this lemma, fixing $\eps\in (0,1/2)$ and taking $L$ as in Part
(1) (we can obviously assume $L\geq L_0$, since otherwise $z$ is
bounded) we derive from Proposition~\ref{pr-equi} that
$$
\sum_{f\in \Hf_k^*(N)}{\omega_f \Bigl(
A_L(\theta_f/\pi)-B_L(\theta_f/\pi)
\Bigr)}
=\Delta+O(D^{z\pi(z)}(\tau(N)\log 2N)^2(Nk^{5/6})^{-1})
$$
for some absolute constants $D$, with $\Delta\geq
(1/2-\eps)^{\pi(z)}$. This is then $\gg (1/2-\eps)^{\pi(z)}$, provided
$$
D^{z\pi(z)}(\tau(N)\log 2N)^2(Nk^{5/6})^{-1}\ll \Bigl(\frac{1}{2}-
\eps\Bigr)^{\pi(z)}.
$$
\par
This condition is satisfied for
$$
z\leq c\sqrt{(\log kN)(\log\log kN)}
$$
where $c>0$ is an absolute constant, and this gives
Theorem~\ref{th-lower-bound} when counting with the weight $\omega_f$.
But, using well-known bounds for $\langle f,f\rangle$, we have
$$
\omega_f\ll kN(\log kN)(\log\log 6N)\ll kN(\log kN)^2,
$$
with an absolute implied constant. Hence, for $z=c\sqrt{(\log
  kN)(\log\log kN)}$, we get
\begin{align*}
\frac{1}{|\Hf_k^*(N)|} |\{f\in \Hf_k^*(N)\,\mid\, \lambda_f(p)\geq 0\text{ for } p\leq z,\ p\nmid N\}| &\gg
\frac{1}{(\log kN)^2}\Bigl(\frac{1}{2}- \eps\Bigr)^{\pi(z)}\\
&\gg \Bigl(\frac{1}{2}-2\eps\Bigr)^{\pi(z)}
\end{align*}
if $kN$ is large enough, and so we obtain Theorem~\ref{th-lower-bound}
as stated.

\begin{proof}[Proof of Lemma~\ref{lm-polys}]
  We must now refer to the specific construction in~\cite{BMV00}. We
  start with $A_L(x)$: we have the product formula
$$
A_L(x)=\prod_{p\mid P}{\alpha_L(x_p)},
$$
where $\alpha_L$ is a trigonometric polynomial in one variable of
degree $\leq L$, i.e., of the type
$$
\alpha_L(x)=\sum_{|\ell|\leq
    L}{\hat{\alpha}_L(\ell)e(\ell x)},
$$
with $\hat{\alpha}_L(0)=1/2$ (see~\cite[(2.2), Lemma 5,
(2.17)]{BMV00}). In particular, the constant term (in the Chebychev
expansion) for $A_L(\theta/\pi)$ is given by
$$
\Bigl(\int_0^{\pi}{\alpha_L(\theta/\pi)\d \mu_{ST}}\Bigr)^{\omega},
$$
and we will bound it below. For the moment, we observe further that,
from~\cite[Lemma 5]{BMV00}, we know that $0\leq \alpha_L(x)\leq 1$ for
all $x\in [0,1]$, and so we can simply bound all the coefficients in
the Chebychev expansion, using the Cauchy-Schwarz inequality and
orthonormality:
\begin{align*}
\Bigl|\int_0^{\pi}{\alpha_L(\theta/\pi)X_n(\theta)\d \mu_{ST}}\Bigr|^2
&\leq \int_{0}^{\pi}
{|\alpha_L(\theta/\pi)|^2\d \mu_{ST}}
\times
\int_{0}^{\pi}
{|X_n(\theta)|^2\d \mu_{ST}}\\
&\leq \int_{0}^{\pi}
{\d \mu_{ST}}
\times
\int_{0}^{\pi}
{|X_n(\theta)|^2\d \mu_{ST}}=1.
\end{align*}
\par
It is also clear using the definition of $X_n(\theta)$ that the
$n$-th coefficient is zero as soon as $n+2>2L$. 
\par
We now come to $B_L(x)$, which is a sum of $\omega$ product functions,
as already indicated: we have
$$
B_L(x)= \sum_{p\mid P}{\beta_L(x_p)\prod_{\stacksum{q\mid
      P}{q\not=p}}{ \alpha_L(x_q) }},
$$
where $\beta_L(x)$ is another trigonometric polynomial of degree $L$,
given explicitly by
\begin{align*}
\beta_L(x)&=
\frac{1}{2L+2}\Bigl(
\sum_{|\ell|\leq L}{
\Bigl(1-\frac{|\ell|}{L+1}\Bigr)e(\ell x)
}
+\sum_{|\ell|\leq L}{
\Bigl(1-\frac{|\ell|}{L+1}\Bigr)e(\ell (x-1/2))
}
\Bigr)\\
&=\frac{1}{2L+2}\Bigl(2+2
\sum_{1\leq \ell\leq L}{
\Bigl(1-\frac{\ell}{L+1}\Bigr)(1+(-1)^{\ell})\cos(2\pi \ell x)
}\Bigr),
\end{align*}
(see~\cite[p. 342, (2.3), p. 339]{BMV00}). 
\par
We now see immediately that Part (3) of the lemma is valid, and
moreover, we see that $|\beta_L(x)|\leq 1$, so the same Cauchy-Schwarz
argument already used for $\alpha_L$ implies that Part (2) holds.
\par
To conclude, we look at the constant term in the Chebychev expansion
for $B_L$, which is given by
$$
\omega
\Bigl(\int_0^{\pi}{\alpha_L(\theta/\pi)\d \mu_{ST}}\Bigr)^{\omega-1}
\int_{0}^{\pi}{\beta_L(\theta/\pi)\d \mu_{ST}}.
$$
\par
Using the expression
$$
\beta_L(\theta/\pi)=\frac{1}{2L+2}\Bigl(2+2
\sum_{1\leq \ell\leq L}{
\Bigl(1-\frac{|\ell|}{L+1}\Bigr)(1+(-1)^{\ell})\cos(2\ell \theta)
}\Bigr),
$$
where the second term doesn't contribute after integrating against
$\sin^2\theta=(1-\cos 2\theta)/2$ (the term with $\ell=1$ is zero), we
get the formula
$$
\Delta=\Bigl(\int_0^{\pi}{\alpha_L(\theta/\pi)\d \mu_{ST}}\Bigr)^{\omega-1}
\Bigl(\int_0^{\pi}{\alpha_L(\theta/\pi)\d \mu_{ST}}
-\frac{\omega}{L+1}\Bigr)
$$
for the contribution of $A_L(x)-B_L(x)$.
\par
Now we come back to a lower bound for the constant term for
$\alpha_L$. The point is that, as $L\ra +\infty$, $\alpha_L$ converges
in $L^2([0,1])$ to the characteristic function $\chi$ of $[0,1/2]$:
from~\cite[(2.6)]{BMV00}, and the definition of $\alpha_L$, we get
$$
|\chi(x)-\alpha_L(x)|\leq \beta_L(x),\quad\quad 0\leq x\leq 1,
$$
and from the Fourier expansion of $\beta_L$ we have
$$
\|\beta_L\|_{L^2}^2\leq \frac{1}{(2L+2)^2}\times (4L+4)\ra 0.
$$ 
\par
Hence, we know that
$$
\frac{2}{\pi}\int_0^{\pi}{\alpha_L(\theta/\pi)\sin^2\theta \d \theta}
\ra \int_0^{\pi}{\chi(\theta/\pi)\d \mu_{ST}}=1/2.
$$
\par
For given $\eps\in (0,1/2)$, the integral is $\geq (1/2-\eps/2)$ if
$L\geq L_0$, for some constant $L_0$. Then, if $L+1\geq
2\eps^{-1}\omega$, we derive
$$
\Delta\geq \Bigl(\frac{1}{2}-\frac{\eps}{2}\Bigr)^{\omega-1}
\Bigl(\frac{1}{2}-\eps\Bigr)\geq 
\Bigl(\frac{1}{2}-\eps\Bigr)^{\omega},
$$
which gives Part (1) of the lemma.
\end{proof}

\section{Proof of Theorem \ref{th-1}}
\label{sec-recognition}

The simple idea of the proof of Theorem~\ref{th-1} is that the
assumption translates to $\lambda_{f_1}(p)\lambda_{f_2}(p)\geq 0$ for
all primes $p$ (with few exceptions). However, it is well-known from
Rankin-Selberg theory that if $f_1\not=f_2$, we have
\begin{equation}\label{eq-rs}
\sum_{p}{\frac{\lambda_{f_1}(p)\lambda_{f_2}(p)}{p^{\sigma}}}=O(1)
\qquad
(\sigma\to 1+)
\end{equation}
(see, e.g,~\cite[\S 5.12]{IK04} for a survey and references; the
underlying fact about automorphic forms is due to M\oe glin and
Waldspurger). Thus we only need to find a lower bound for the
left-hand side (which is a sum of non-negative terms) which is
unbounded as $\sigma$ tends to $1+$. Since Rankin-Selberg theory also
gives
\begin{equation}\label{eq-rs2}
\sum_{p}{\frac{\lambda_{f_1}(p)^2}{p^{\sigma}}}\sim -\log(\sigma-1)
\qquad
(\sigma\to 1+),
\end{equation}
the only difficulty is that one might fear that the coefficients of
$f_1$ and $f_2$ are such that whenever $\lambda_{f_1}(p)$ is not
small, the value of $\lambda_{f_2}(p)$ \emph{is} very
small.\footnote{\ See the remark after the proof for an example of
  which potential situations must be excluded.}  In other words, we
must show that the smaller order of magnitude of~(\ref{eq-rs})
compared with~(\ref{eq-rs2}) is not due to the small size of the
summands, but to sign compensations.  For this we use the following
trick which exploits the little partial information known towards the
pair Sato-Tate conjecture.
\par
Assume first that $f_1$ and $f_2$ are non-CM cusp forms, and that
neither is a quadratic twist of the other (in particular,
$f_1\not=f_2$). By Ramakrishnan's Theorem (\cite[Th. M, \S
3]{ramakrishnan}), there exists a \emph{cuspidal} automorphic
representation $\pi$ on $GL(4)/\Q$ such that
$$
L(\pi,s)=L(f_1\times f_2,s),
$$
and consequently, by Rankin-Selberg theory on $GL(4)\times GL(4)$ now
(the fact that $L(\pi\times\bar{\pi},s)$ has a single pole at $s=1$),
we have
\begin{equation}\label{eq-ramak}
\sum_{p}{\frac{(\lambda_{f_1}(p)\lambda_{f_2}(p))^2}{p^{\sigma}}}
= \sum_{p}{\frac{1}{p^{\sigma}}}+O(1),
\qquad
(\sigma\to 1+).
\end{equation}
\par
However, if we denote by $E$ the set of primes $p$ for which
$\lambda_f(p)\lambda_g(p)<0$, we have
\begin{align*}
  \sum_{p}{\frac{(\lambda_{f_1}(p)\lambda_{f_2}(p))^2}{p^{\sigma}}}
  &=\sum_{p\notin
    E}{\frac{(\lambda_{f_1}(p)\lambda_{f_2}(p))^2}{p^{\sigma}}}
  +\sum_{p\in
    E}{\frac{(\lambda_{f_1}(p)\lambda_{f_2}(p))^2}{p^{\sigma}}}\\
  &\leq 4\sum_{p\notin
    E}{\frac{\lambda_{f_1}(p)\lambda_{f_2}(p)}{p^{\sigma}}}
+16\sum_{p\in E}{\frac{1}{p^{\sigma}}}
\end{align*}
by Deligne's bound. Then the first sum can also be written
\begin{align*}
4\sum_{p\notin
    E}{\frac{\lambda_{f_1}(p)\lambda_{f_2}(p)}{p^{\sigma}}}&=
-4\sum_{p\in E}{\frac{\lambda_{f_1}(p)\lambda_{f_2}(p)}{p^{\sigma}}}
+\sum_{p}{\frac{\lambda_{f_1}(p)\lambda_{f_2}(p)}{p^{\sigma}}}\\
&\leq 16 \sum_{p\in E}{\frac{1}{p^{\sigma}}}+O(1)
\end{align*}
using once more Deligne's bound and the assumption $f_1\not=f_2$ to
apply~(\ref{eq-rs}). 
\par
Comparing~(\ref{eq-ramak}) with these two inequalities leads to
$$
\sum_{p\in E}{\frac{1}{p^{\sigma}}}\geq
  \frac{1}{32}\sum_{p}\frac{1}{p^{\sigma}}
+O(1),\quad\quad \sigma\ra 1,
$$
i.e., the set of primes where the signs of $f_1$ and $f_2$ differ has
analytic density $\geq 1/32$.
\par
There remains to consider Part (1) of Theorem~\ref{th-1} when one of
the forms is of CM type (and the exceptional set $E$ now has density
$0$). We will be brief since there are less difficulties here. The
main point is the following well-known result concerning the
distribution of the angles $\theta_f(p)$ for a CM form $f\in
\Hf_k^*(N)$, with $k\geq 2$: there exists a real, non-trivial,
primitive Dirichlet character $\chi_f$ such that $\lambda_f(p)=0$ when
$\chi_f(p)=-1$ (a set of primes $I_f$ of density $1/2$), and for
$p\notin I_f$, the $\theta_f(p)\in [0,\pi]$ for $p\leq x$ become
\emph{uniformly distributed} as $x\ra +\infty$, i.e., we have
$$
\frac{2}{\pi(x)}\sum_{\stacksum{p\notin I_f}
{p\leq x}}{e^{2im\theta_f(p)}}\ra 0,
$$
for all non-zero integers $m\in\Z$ (see, e.g.,~\cite[p. 197]{km},
where this is explained for elliptic curves, with slightly different
notation). In particular, for any $\alpha>0$, the density of the set
of primes where $|\lambda_f(p)|>\alpha$ exists and is equal to
$$
\frac{1}{\pi}\arccos(\alpha/2)
$$
and this density goes to $1/2$ as $\alpha\ra 0$. 
\par
Now assume $f_1$ is a CM form and $f_2$ is not; according to
Lemma~\ref{lm-non-cm-case} below, we find $\alpha>0$ and a set of
primes $P_2$ of analytic density $\delta>1/2$ where
$|\lambda_{f_2}(p)|>\alpha$, and then the set $P_2\cap I_{f_1}$ has
analytic density $>0$, thus for small enough $\alpha'$, it contains a
set $G$ with positive analytic density where
$|\lambda_{f_1}(p)|>\alpha'$. Hence we have
\begin{align*}
    \sum_{p}
  {\frac{\lambda_{f_1}(p)\lambda_{f_2}(p)}{p^{\sigma}}}
  & \geq  \sum_{p\in G}
  \frac{\lambda_{f_1}(p)\lambda_{f_2}(p)}{p^{\sigma}} + o(\log|\sigma-1|^{-1})
  \\
  &\geq \alpha\alpha'\sum_{p\in
    G}{\frac{1}{p^{\sigma}}}+o(\log|\sigma-1|^{-1}),
\quad\text{ as } \sigma\ra 1+,
\end{align*}
which is in fact a contradiction (since $f_1$ can not be equal to
$f_2$).
\par
Finally, assume $f_1$ and $f_2$ are CM forms. Because of independence
of primitive real characters, the union $I_{f_1}\cup I_{f_2}$ has
density \emph{at most} $3/4$ (the complement contains the set of
primes totally split in a Galois extension of $\Q$ of degree at most
$4$). For small enough $\alpha>0$, the complement must contain a set
of primes of positive analytic density where
$|\lambda_{f_1}(p)|>\alpha$, $|\lambda_{f_2}(p)|>\alpha$, and we can
conclude as before that the Rankin-Selberg convolution has a pole at
$s=1$, so that $f_1=f_2$ in that case also.

\begin{remark}
In the first version of this paper, we did not use Ramakrishnan's
theorem, but managed to prove a weaker version of Part (1) of
Theorem~\ref{th-1} using only the Rankin-Selberg properties of $f_1$
and $f_2$ together with the analytic properties of (small) symmetric
square $L$-functions.  We sketch the argument, since this may be of
interest in other contexts.
\par
The basic point is the following lemma, which may be of independent
interest: 

\begin{lemma}\label{lm-non-cm-case}
 Let $N\geq 1$ be an integer, $k\geq 2$ be an even integer and $f\in
 {\rm H}_{k}^*(N)$ a primitive cusp form of level $N$ and weight $k$ which
 is not of CM type.  Then there exists a constant $\alpha>0$ and
 $\delta>\dm$ such that
$$
\sum_{|\lambda_f(p)|>\alpha}{\frac{1}{p^{\sigma}}}
 \geq \delta\sum_{p}{\frac{1}{p^{\sigma}}}+O(1),
$$
for $\sigma>1$. In fact, one can take $\alpha=0.231$ and
$\delta=\dm+\frac{1}{24}$.
\end{lemma}



\begin{proof}
  It is convenient here to work with the Chebychev polynomials $U_n$
  instead of the Chebychev functions $X_n$ considered in the previous
  section: recall that for $n\geq 0$, we have
$$
X_n(\theta)=U_n(2\cos\theta)
$$
where $U_n\in \R[x]$ is a polynomial of degree
$n$. Then~(\ref{eq-cheby}) 
gives $U_n(\lambda_f(p))=\lambda_f(p^n)$ for any $f\in\Hf_k^*(N)$,
$p\nmid N$, and $n\geq 0$.
\par
We then claim that there exists a polynomial
$$
Y=\beta_0+\beta_2U_2+\beta_4U_4+\beta_6U_6\in \R[x]
$$
with the following properties: 
\par
(i) 
$\beta_0>\dm$; 
\par
(ii) 
for some $\alpha>0$ and $x\in [-2,2]$, we have 
\begin{equation}\label{eq-approx}
Y(x)\leq \chi_{A}(x),
\end{equation}
where $A := \{x\in [-2,2]\,\mid\,|x|>\alpha\}$.
\par
Assuming this, we conclude as follows: by (ii),
we have
$$
\sum_{|\lambda_f(p)|>\alpha}{\frac{1}{p^{\sigma}}}
\geq \sum_{p\nmid N}{\frac{Y(\lambda_f(p))}{p^{\sigma}}}
= \beta_0\sum_{p\nmid N}{\frac{1}{p^{\sigma}}}
+\sum_{1\le i\le 3} {\beta_{2i}\sum_{p\nmid N}{\frac{U_{2i}(\lambda_f(p))}{p^{\sigma}}}}.
$$
By the holomorphy and non-vanishing at $s=1$ of the second, fourth and
sixth symmetric power $L$-functions (see~\cite[Th. 3.3.7,
Prop. 4.3]{kim-shahidi} for the last two, noting that non-CM forms are
not dihedral, and~\cite{shahidi} for a survey concerning those
$L$-functions), since $U_n(\lambda_f(p))$ is exactly the $p$-th
coefficient of the $n$-th symmetric power for $p\nmid N$, standard
analytic arguments show that
$$
\sum_{p\nmid N}{\frac{U_{2i}(\lambda_f(p))}{p^{\sigma}}}=O(1)
$$
for $\sigma\geq 1$ and $i=1,2,3$. Hence the result follows with
$\delta=\beta_0>\dm$.
\par
Now to check the claim, and verify the values of $\alpha$ and
$\delta$, we just exhibit a suitable polynomial, namely
$$
Y
 = \textstyle \frac{1}{2}+\frac{1}{24}+\frac{1}{4}U_2-
\frac{1}{4}U_4+\frac{136}{1000}U_6
 = \textstyle \frac{17}{125}x^6 - \frac{93}{100}x^4 + \frac{227}{125}x^2 -
\frac{283}{3000},
$$
since
\begin{equation}\label{eq-first-xs}
\begin{cases}
U_0 =1, \quad
U_1=x,\quad 
U_2 =x^2-1,
\quad
U_4
= x^4-3x^2+1, \\
U_5 
= x^5-4x^3+3x,\quad
U_6=x^6-5x^4+6x^2-1.
\end{cases}
\end{equation}
\par
This polynomial is even, and its graph on $[-2,2]$ is in Figure~1.
\par
\begin{figure}[ht]
\centering
\includegraphics[width=4in]{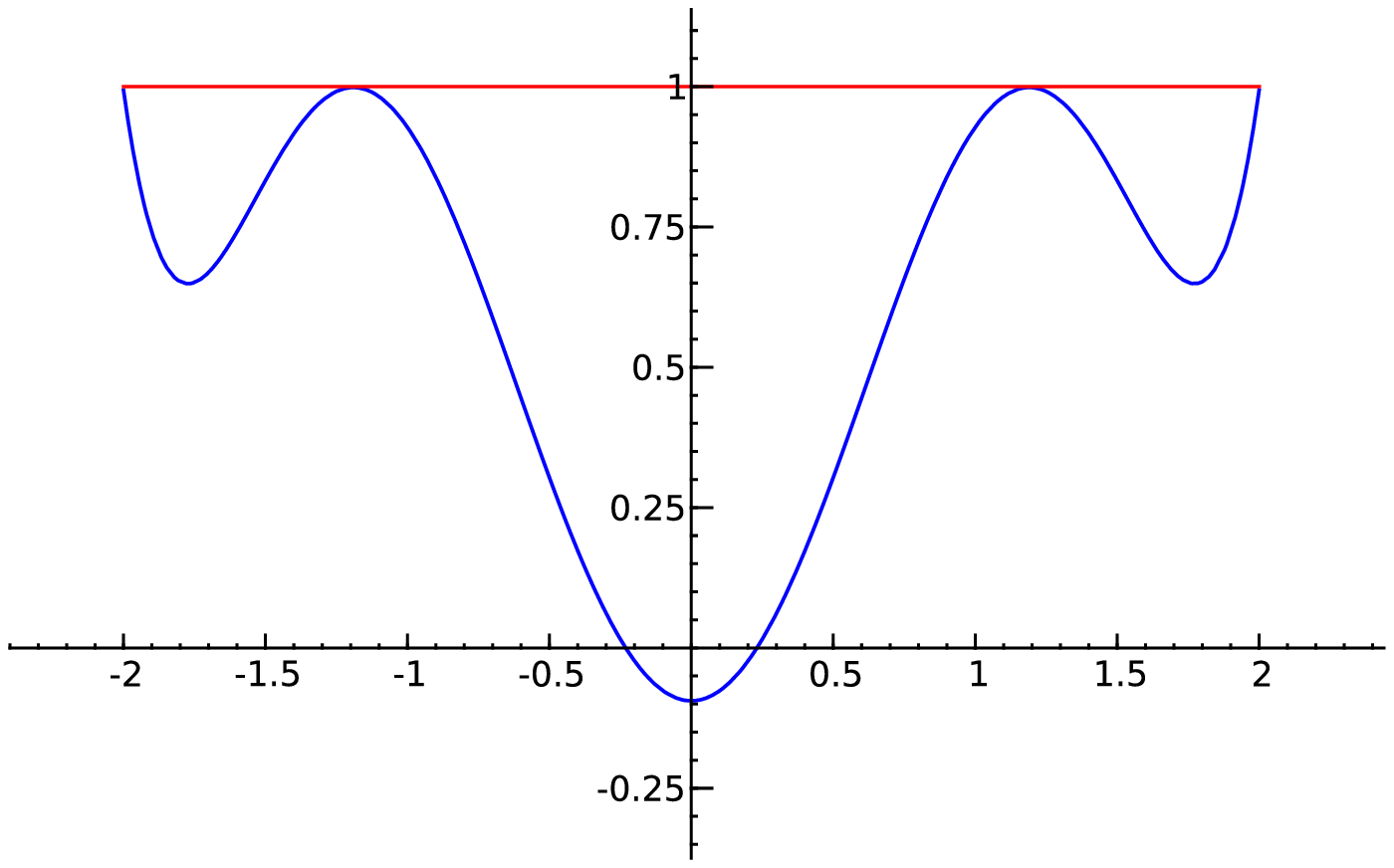}
\caption{}
\end{figure}
The value of $\alpha$ is an approximation (from below) to the real
root
$$
\alpha_0=0.23107202470801418176315245050693402580\ldots
$$
of $Y$ in $[0,2]$; the maximum value of $Y$ on $[0,2]$ is very close
to $1$.
\end{proof}

The upshot of this lemma is that, given $f_1$ and $f_2$ (not of CM
type), there exists a set of primes of analytic density $>0$ for which
both $|\lambda_{f_1}(p)|$ and $|\lambda_{f_2}(p)|$ have a positive
lower bound. Then the sum of $\lambda_{f_1}(p)\lambda_{f_2}(p)$ over
this set can not be small, and this leads to an upper bound for the
density of the ``exceptional set''. However, the actual value from the
above lemma is much smaller than what Theorem~\ref{th-1} uses (it is
about $1/1000$).

Another interesting point of this method is that using the sixth
symmetric power (and thus the deep results of Kim and Shahidi) is
necessary for Lemma~\ref{lm-non-cm-case}. For this, note that the
sequences $\{x_p\}_{p\,{\rm primes}}$ and $\{y_p\}_{p\,{\rm primes}}$
defined by $x_2=y_2=0$ and for primes $p\geq 3$ by
\begin{align*}
  x_p & = \begin{cases} 0 & \text{if $p\equiv 3\mods{4}$},
    \\
    (-1)^{(p-1)/4}\sqrt{2} & \text{if $p\equiv 1\mods{4}$},
\end{cases}
\\
y_p & =\begin{cases} (-1)^{(p-3)/4}\sqrt{2} & \text{if $p\equiv
    3\mods{4}$},
  \\
  0 & \text{if $p\equiv 1\mods{4}$},
\end{cases}
\end{align*}
have the ``right'' moments of order $1$ to $5$ for being Sato-Tate
distributed,\footnote{\ The sixth moment fails: it is $4$ instead of
  $5$ for the Sato-Tate distribution.} i.e., we have
$$
\sum_{p} \frac{X_k(x_p)}{p^{\sigma}}
= O(1)
\qquad{\rm and}\qquad
\sum_{p} \frac{X_k(y_p)}{p^{\sigma}}
=O(1)
$$
for $\sigma>1$ and $1\leq k\leq 5$,
and yet $x_py_p\geq 0$ for all $p$, in fact $x_py_p=0$, so that we
most certainly have
$$
\sum_{p} \frac{x_py_p}{p^{\sigma}}=O(1)
\qquad(\sigma>1).
$$
\end{remark}

\begin{remark}
  As a final remark, one can think of other ways (than looking at
  signs) of reducing Fourier coefficients of modular forms to a fixed
  finite set: the most obvious, at least if $f$ has integral
  coefficients $\lambda_f(n)n^{(k-1)/2}$, is to look at the
  coefficients modulo some fixed prime number $\ell$. However, the
  situation there can be drastically different: for instance, for all
  (infinitely many) elliptic curves with full rational $2$-torsion,
  given for instance by equations
$$
y^2=(x-a)(x-b)(x-c)
$$
with $a$, $b$, $c$ distinct integers, the reduction modulo $2$ of the
odd prime coefficients of the corresponding $L$-function (or modular
form) is the same!
\end{remark}

\end{document}